\documentclass[smallextended,referee,envcountsect,]{svjour3}

\smartqed
\usepackage{graphicx}
\usepackage{amsmath}
\usepackage{amsfonts}
\usepackage{amssymb}
\usepackage{bm}

\usepackage{tikz}
\usepackage{algorithm}
\usepackage{algpseudocode}

\usepackage{booktabs}%
\usepackage{multirow}%

\newtheorem{assumption}{Assumption}
\journalname{JOTA}

\begin{document}

\title{Computing the Minimum-Time Interception of a Moving Target}

\author{Maksim Buzikov}

\institute{Maksim Buzikov \at
             V.A.~Trapeznikov Institute of Control Sciences of Russian Academy of Sciences \\
             Moscow, Russia\\
             me.buzikov@physics.msu.ru
}

\date{Received: date / Accepted: date}

\maketitle

\begin{abstract}
    In this study, we propose an algorithmic framework for solving a class of optimal control problems. This class is associated with the minimum-time interception of moving target problems, where a plant with a given state equation must approach a moving target whose trajectory is known a priori. Our framework employs an analytical description of the distance from an arbitrary point to the reachable set of the plant. The proposed algorithm is always convergent and cannot be improved without losing the guarantee of its convergence to the correct solution for arbitrary Lipschitz continuous trajectories of the moving target. In practice, it is difficult to obtain an analytical description of the distance to the reachable set for the sophisticated state equation of the plant. Nevertheless, it was shown that the distance can be obtained for some widely used models, such as the Dubins car, in an explicit form. Finally, we illustrate the generality and effectiveness of the proposed framework for simple motions and the Dubins model.
\end{abstract}

\section{Introduction}\label{introduction}

The problem of intercepting a moving target (PIMT) arises in both civilian and military problem areas. Most modern guidance laws for real problems are derived using linear quadratic optimal control theory to obtain explicit analytical feedback solutions~\cite{Palumbo2010-sd}. This approach requires linear approximation of the problem. Another way to solve a PIMT is to simplify the state equation such that the corresponding optimal control problem can have a closed-form solution. The solution to the relaxed problem produces a reference path for a real plant. This requires the use of path-following control strategies~\cite{Rubi2020-ts}. Relaxation of the problem requires the selection of a state equation that simultaneously satisfies several requirements. First, in practice, the reference path provided by the simplified model should be close to a feasible path. Secondly, path computation should be a time-rapid and memory-efficient process for onboard computers. The main motivation of this study was to investigate a compromise method for PIMT.

In the literature, numerous models have been studied analytically. For example, isotropic rocket~\cite{Isaacs1965-nd,Lewin1989-ly,Akulenko1998-mr,Akulenko2001-xw,Akulenko2018-ck,Venkatraman2006-nc,Bakolas2014-ze}, Dubins' model~\cite{Markov1889-iu,Dubins1957-hl,Pecsvaradi1972-xl,Shkel2001-gf,Bui1993-qt,Patsko2003-nu,Kaya2017-zp,Coates2017-fl} (including asymmetric cases~\cite{Bakolas2011-jk,Patsko2018-ef}), Dubins airplane~\cite{Chitsaz2007-pq,McLain2014-za,Vana2020-kz}, and Reeds-Shepp model~\cite{Reeds1990-nc,Sussmann1991-wt,Boissonnat1994-xv,Soueres1996-vk}. Such models, in comparison with linearly approximated models, can consider the maneuverability of a real plant more accurately. At the same time, these models retain the possibility of analytically solving the problem of computing optimal trajectories.

In this study, we considered a class of problems involving the minimum time interception of a target that moves along a known and predetermined trajectory. We assume that the state equation of the plant is so simple that we can analytically describe the distance from an arbitrary point to the reachable set of the plant. The only restriction imposed on the trajectory of a moving target is that it is a Lipschitz continuous function of the time. From a practical perspective, this implies that target coordinates vary at a limited rate. Specific cases of such a problem have been reported in several studies. In~\cite{Bakolas2014-ze}, the problem of intercepting a moving target using an isotropic rocket is explored. In a series of studies~\cite{Clements1990-fp,Looker2008-yx,Zhang2014-fp,Meyer2015-st,Zheng2021-wu}, the same problems were investigated for the interception of a uniformly moving target using a Dubins car. A more general problem of intercepting by a Dubins car for an arbitrary continuous target trajectory was considered in~\cite{Zheng2021-wu,Buzikov2021-md}. The lateral interception of a moving target by a Dubins car was explored in~\cite{Techy2009-jv,Bakolas2013-zs,Gopalan2017-wy,Mittal2020-gh,Zheng2022-bh} for a uniformly moving target, in~\cite{Manyam2020-ym,Manyam2022-gh} for a target moving on a circle or a racetrack path, and in~\cite{McNeely2007-hw,Buzikov2022-ff} for an arbitrary continuous trajectory of the target. The main practical contribution of this study is the always convergent algorithm for calculating the minimum time interception of a moving target by a Dubins car. In comparison with~\cite{Zhang2014-fp,Zheng2021-wu}, the proposed algorithm works not only with rectilinear uniform movement of the target. It also calculates the minimum time required for arbitrary Lipschitz continuous trajectories (including rectilinear uniform trajectories).

The analytical description of a reachable set in a simple motion model is trivial. In more complex cases, such as Dubins' model, a series of analytical results is obtained. For example, an analytical description of a planar reachable set of a Dubins car can be found in~\cite{Cockayne1975-pu,Boissonnat1994-pd,Fedotov2011-gr,Cacace2020-pu,Buzikov2021-md}. Papers~\cite{Patsko2003-nu,Patsko2018-ef,Patsko2022-pl,Buzikov2022-ff} have been devoted to the study of a three-dimensional reachable set of a Dubins car. In this study, we assume that the distance from an arbitrary point to the reachable set of the plant is known or can be calculated effectively. If the distance to the target position is calculated, the minimum time interception is such that the distance is not greater than the desired value of the capture radius. The main idea of designing an always convergent algorithm to obtain the optimal interception time was borrowed from~\cite{Chernousko1968-nd,Sukharev1976-pb,Abaffy2013-rf,Galantai2015-xx}. To specify this idea, we additionally used the properties of distance to the reachable set. The novelty of this study in comparison with these works lies in the use of additional information regarding the description of the reachable set to increase the step of the fixed-point iteration algorithm.

The remainder of this paper is organized as follows. In Section~\ref{sec:problem_formulation}, we formally describe the minimum time interception of a moving target problem. In Section~\ref{sec:reachable_sets}, we investigate some properties of reachable sets and universal lower estimators of interception time. In addition, we proposed an iterative algorithm based on universal lower estimators. In Section~\ref{sec:examples}, we demonstrate the proposed algorithmic framework using two examples of plants: a simple motion and Dubins car. Additionally, we present numerical experiments for these two examples. Finally, Section~\ref{sec:conclusion} concludes the paper with a discussion and brief description of future work.

\section{Problem Formulation}\label{sec:problem_formulation}

Let $t \in \mathbb{R}^+_0$ denote a time moment and $\boldsymbol{x} = \mathrm{stack}(\boldsymbol{y}, \boldsymbol{z})$ denote the state vector function. Here, $\boldsymbol{x}(t) \in \mathcal{X} = \mathcal{Y} \oplus \mathcal{Z}$. $\mathcal{X}$ is a state space. $\mathcal{Y}$ is a finite-dimensional normed space with norm $\lVert \cdot \rVert: \mathcal{Y} \to \mathbb{R}^+_0$. We assume that $\boldsymbol{y}(t) \in \mathcal{Y}$ corresponds to coordinates that are important for interception. The remaining coordinates $\boldsymbol{z}(t) \in \mathcal{Z}$ are arbitrary when the target is intercepted. The state equation (dynamic constraints) of the plant is as follows:
\begin{equation}
    \dot{\boldsymbol{y}} = 
    \boldsymbol{f}(\boldsymbol{x}, \boldsymbol{u}), \quad \dot{\boldsymbol{z}} = \boldsymbol{g}(\boldsymbol{x}, \boldsymbol{u}).
\end{equation}
Here, $\boldsymbol{u}(t) \in \mathcal{U}$ is a control input and $\mathcal{U}$ is a compact restraint set of admissible values of the control inputs. The class of all measurable admissible control inputs is denoted as $\mathcal{A}$ ($\boldsymbol{u} \in \mathcal{A}$). Throughout this paper, we assume that the control input and state equation are subject to the requirements of Theorem 2 from~\cite[pp. 242--244]{Lee1986-vf}. In addition, we assumed that the velocity of the plant is uniformly restricted. Without loss of generality, we claim $\lVert \boldsymbol{f}(\mathbf{x}, \mathbf{u}) \rVert \leq 1$ for all\footnote{We will use italic letters $\boldsymbol{x}$, $\boldsymbol{y}$, $\boldsymbol{u}$ ... for functions and roman letters $\mathbf{x}$, $\mathbf{y}$, $\mathbf{u}$ ... for points of the corresponding spaces.} $\mathbf{x} \in \mathcal{X}$ and $\mathbf{u} \in \mathcal{U}$.

\begin{figure}
    \centering
    \begin{tikzpicture}[scale=1.5]
        \draw[->] (0, -0.25) -- (0, 3) node[above]{$\mathcal{Y}$};
        \draw[->] (0, 0) -- (5, 0);
        \filldraw (0, 0) circle (1pt) node[left]{$0$};
        \draw[-] (0, 0) .. controls (2, 2) and (3, 1.5) .. (4, 2);
        \draw[dashed,yshift=-6pt] (0, 0) .. controls (2, 2) and (3, 1.5) .. (4, 2);
        \draw[dashed,yshift=6pt] (0, 0) .. controls (2, 2) and (3, 1.5) .. (4, 2);
        \node at (1.95, 1.5) {$\boldsymbol{y}(\cdot)$};
        \draw[densely dotted] (0, 2) .. controls (2, 2.5) and (3, 2.25) .. node[above]{$\boldsymbol{y}_T(\cdot)$} (4, 2.25);
        \filldraw[black] (4, 2.25) circle (1pt) node[right]{$\lVert \boldsymbol{y}(t) - \boldsymbol{y}_T(t) \rVert \leq \ell$};
        \filldraw[black] (4, 0) circle (1pt) node[below]{$t$};
        \draw[dotted] (4, 0) -- (4, 2.25);
    \end{tikzpicture}
    \caption{Trajectory of interception (solid line) of the moving target (dotted line). The dashed lines around the trajectory of interception define the capture set (radius $\ell$). The target enters the capture set at $t$. Space $\mathcal{Y}$ is depicted as a line for simplicity}
    \label{fig:problem}
\end{figure}
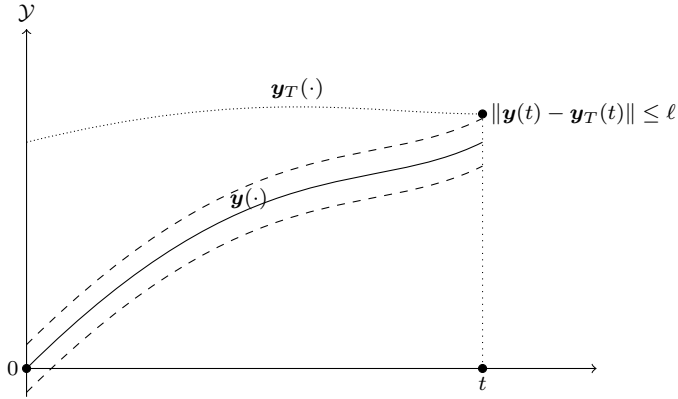

Let $\boldsymbol{y}_T \in \mathrm{Lip}_v(\mathbb{R}^+_0, \mathcal{Y})$\footnote{$\mathrm{Lip}_v(\mathcal{F}, \mathcal{G})$ denotes a set of $v$-Lipschitz continuous functions. If $\boldsymbol{y} \in \mathrm{Lip}_v(\mathcal{F}, \mathcal{G})$, then $\boldsymbol{y}: \mathcal{F} \to \mathcal{G}$, and $\lVert \boldsymbol{y}(\boldsymbol{x}_2) - \boldsymbol{y}(\boldsymbol{x}_1) \rVert_\mathcal{G} \leq v \lVert \boldsymbol{x}_2 - \boldsymbol{x}_1 \rVert_\mathcal{F}$ for all $\boldsymbol{x}_1, \boldsymbol{x}_2 \in \mathcal{F}$.} describe the trajectory of the target (Fig.~\ref{fig:problem}). The constant $v \in \mathbb{R}^+_0$ denotes the maximum speed of the target. We assumed that the trajectory of the target is known a priori. However, the algorithm for obtaining the solution should work with any trajectory $\boldsymbol{y}_T \in \mathrm{Lip}_v(\mathbb{R}^+_0, \mathcal{Y})$ and use only knowledge about the maximal speed of the target.

Let $\boldsymbol{x}(0) = \mathbf{x}_0 \in \mathcal{X}$. Without loss of generality, we set $\boldsymbol{y}(0) = \boldsymbol{0} \in \mathcal{Y}$. Unless stated otherwise, we assume that $\boldsymbol{x} = \mathrm{stack}(\boldsymbol{y}, \boldsymbol{z})$ is an absolutely continuous solution of the state equation with a control input $\boldsymbol{u} \in \mathcal{A}$ and initial conditions $\boldsymbol{x}(0) = \mathbf{x}_0$. We define the problem of the minimum-time interception of a moving target by finding the minimum value of the following functional:
\begin{equation}
    J[\boldsymbol{u}; \boldsymbol{y}_T] \overset{\mathrm{def}}{=} \min\left\{t \in \mathbb{R}^+_0:\: \lVert \boldsymbol{y}(t) - \boldsymbol{y}_T(t) \rVert \leq \ell\right\} \to \inf_{\boldsymbol{u} \in \mathcal{A}}.
\end{equation}
Here, $\ell \in \mathbb{R}^+_0$ denotes the capture radius. If $\ell = 0$, then the interception of the target is a coincidence of the plant and target $\mathbf{y}$-coordinates. Throughout this paper, it is convenient to define $\min \varnothing = +\infty$ for the uniformity of the notation.

\section{Interception by Reachable Sets}\label{sec:reachable_sets}

\begin{figure}
    \centering
    \begin{tikzpicture}
        \draw[->] (0, 0) -- (0, 3) node[above]{$\mathcal{Y}$};
        \draw[->] (0, 0) -- (-1, -1.5) node[below]{$\mathcal{Z}$};
        \draw[->] (0, 0) -- (4, 0);
        \filldraw (0, 0) circle (1pt);
        \filldraw (-1/3, -0.5) circle (1pt) node[left]{$\mathbf{x}_0$};
        \draw[-] (-1/3, -0.5) .. controls (1, 0) and (2.5, 1.2) .. (3, 1.5);
        \draw[-] (-1/3, -0.5) .. controls (1, 0) and (2.5, 0.3) .. (3, 0.5);
        \draw[dotted] (-1/3, -0.5) .. controls (1, 0) and (2.5, 1.0) .. (3, 1.3);
        \draw[dotted] (-1/3, -0.5) .. controls (1, 0) and (2.5, 0.8) .. (3, 1.0);
        \draw[dotted] (-1/3, -0.5) .. controls (1, 0) and (2.5, 0.6) .. (3, 0.7);
        \draw[dashed] (3, 1.0) ellipse (0.1 and 0.5);
        \filldraw (3.3, 0) circle (1pt) node[above right]{$t$};
        \draw[dashed] (3.3, 0) -- (2.3, -1.5);
        \draw[dashed] (3.3, 0) -- (3.3, 3);
        \draw[densely dotted] (3, 0.5) -- (3.3, 0.95);
        \draw[densely dotted] (3, 1.5) -- (3.3, 1.95);
        \draw[-] (3.3, 0.95) -- (3.3, 1.95);
        \draw[densely dotted] (0, 0.95) -- (3.3, 0.95);
        \draw[densely dotted] (0, 1.95) -- (3.3, 1.95);
        \draw[line width=1] (0, 0.95) -- node[left]{$\mathcal{R}(t)$} (0, 1.95);
        \filldraw (0, 0.95) circle (1pt);
        \filldraw (0, 1.95) circle (1pt);
    \end{tikzpicture}
    \caption{$\mathcal{R}(t)$ (bold line) is a projection of the reachable set (solid line) at time $t$ (dashed ellipse) on $\mathcal{Y}$ space}
    \label{fig:R}
\end{figure}
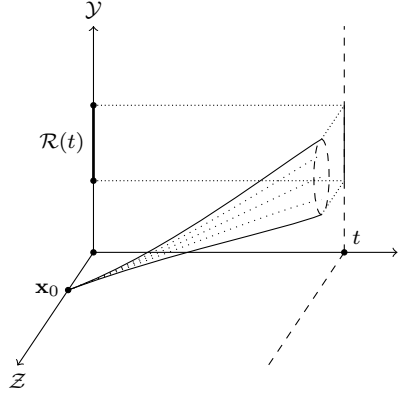

A reachable set at fixed time $t$ is a set of all states in the state space that can be reached at time $t$ from a given initial state by employing an admissible control input~\cite{Patsko2003-nu}. Further, we use the projection of the reachable set on $\mathcal{Y}$ (Fig.~\ref{fig:R}):
\begin{equation}
    \mathcal{R}(t) \overset{\mathrm{def}}{=} \left\{\int\limits_{[0, t]} \boldsymbol{f}(\boldsymbol{x}(\tau), \boldsymbol{u}(\tau)) \mathrm{d}\tau:\: \boldsymbol{u} \in \mathcal{A}\right\}.
\end{equation}

Note that $\mathcal{R}: \mathbb{R}^+_0 \to 2^\mathcal{Y}$ is a multivalued mapping, and $\mathcal{R}(t) \subset \mathcal{Y}$. If $\mathcal{Z} = \varnothing$, $\mathcal{R}(t) \subset \mathcal{Y} = \mathcal{X}$ is a reachable set in the classical sense of the word. According to the general results of optimal control theory, $\mathcal{R}(t)$ is a compact set~\cite[Theorem 2, pp.~242--244]{Lee1986-vf}.

Further, we use the distance function from $\mathcal{R}(t)$ to $\mathbf{y} \in \mathcal{Y}$:
\begin{equation}
    \rho(t, \mathbf{y}) \overset{\mathrm{def}}{=} \min_{\boldsymbol{\eta} \in \mathcal{R}(t)} \lVert \boldsymbol{\eta} - \mathbf{y} \rVert.
\end{equation}
We assume that this function is available in an explicit analytical form for a given plant or that it can be calculated rapidly. This requirement is sufficiently strong, but we demonstrate later that the distance function can be obtained explicitly for the simple motion model and the Dubins model. Thus, the proposed algorithmic framework can be used effectively at least for these motion models.

\begin{lemma}\label{lem:rho_t_lip}
    $\rho(\cdot,\mathbf{y}) \in \mathrm{Lip}_1(\mathbb{R}^+_0, \mathbb{R}^+_0)$ for any $\mathbf{y} \in \mathcal{Y}$.
\end{lemma}
\begin{proof}
    Let $t_1, t_2 \in \mathbb{R}^+_0$. For all $\mathbf{y}_1 \in \mathcal{R}(t_1)$ and $\mathbf{y}_2 \in \mathcal{R}(t_2)$,
    \begin{equation}\label{eq:int1}
        \rho(t_1, \mathbf{y}) = \min_{\boldsymbol{\eta} \in \mathcal{R}(t_1)} \lVert \boldsymbol{\eta} - \mathbf{y} \rVert \leq \lVert \mathbf{y}_1 - \mathbf{y} \rVert \leq \lVert \mathbf{y}_2 - \mathbf{y}_1 \rVert + \lVert \mathbf{y} - \mathbf{y}_2 \rVert.
    \end{equation}
    Using reachability properties, for all $\mathbf{y}_2 \in \mathcal{R}(t_2)$, there exists $\boldsymbol{u} \in \mathcal{A}$ such that 
    \begin{equation}
        \tilde{\mathbf{y}} \overset{\mathrm{def}}{=} \mathbf{y}_2 - \int\limits_{[t_1, t_2]} \boldsymbol{f}(\boldsymbol{x}(t), \boldsymbol{u}(t)) \mathrm{d}t \in \mathcal{R}(t_1).
    \end{equation}
    Since~\eqref{eq:int1} holds for all $\mathbf{y}_1 \in \mathcal{R}(t_1)$, we can substitute $\mathbf{y}_1 = \tilde{\mathbf{y}}$:
    \begin{equation}
        \rho(t_1, \mathbf{y}) \leq \left\lVert\: \int\limits_{[t_1, t_2]} \boldsymbol{f}(\boldsymbol{x}(t), \boldsymbol{u}(t)) \mathrm{d}t \right\rVert + \lVert \mathbf{y} - \mathbf{y}_2 \rVert \\
        \leq \lvert t_2 - t_1 \rvert + \lVert \mathbf{y} - \mathbf{y}_2 \rVert.
    \end{equation}
    The last inequality follows from $\lVert \boldsymbol{f}(\boldsymbol{x}(t), \boldsymbol{u}(t)) \rVert \leq 1$. Since the inequality holds for all $\mathbf{y}_2 \in \mathcal{R}(t_2)$, we can use $\mathbf{y}_2 \in \mathrm{arg}\min_{\boldsymbol{\eta} \in \mathcal{R}(t_2)} \lVert \boldsymbol{\eta} - \mathbf{y} \rVert$. This yields $\rho(t_1, \mathbf{y}) - \rho(t_2, \mathbf{y}) \leq \lvert t_2 - t_1 \rvert$. Similarly, we obtain $\rho(t_2, \mathbf{y}) - \rho(t_1, \mathbf{y}) \leq \lvert t_2 - t_1 \rvert$.\qed
\end{proof}
\begin{lemma}\label{lem:rho_y_lip}
    $\rho(t,\cdot) \in \mathrm{Lip}_1(\mathcal{Y}, \mathbb{R}^+_0)$ for any $t \in \mathbb{R}^+_0$.
\end{lemma}
\begin{proof}
    For all $\mathbf{y}_1, \mathbf{y}_2 \in \mathcal{Y}$,
    \begin{equation}
        \rho(t, \mathbf{y}_1) = \min_{\boldsymbol{\eta} \in \mathcal{R}(t)} \lVert \boldsymbol{\eta} - \mathbf{y}_1 \rVert \leq \min_{\boldsymbol{\eta} \in \mathcal{R}(t)} \lVert \boldsymbol{\eta} - \mathbf{y}_2 \rVert + \lVert \mathbf{y}_2 - \mathbf{y}_1 \rVert = \rho(t, \mathbf{y}_2) + \lVert \mathbf{y}_2 - \mathbf{y}_1 \rVert.
    \end{equation}
    Therefore, $\rho(t, \mathbf{y}_1) - \rho(t, \mathbf{y}_2) \leq \lVert \mathbf{y}_2 - \mathbf{y}_1 \rVert$. In the same way, we can obtain $\rho(t, \mathbf{y}_2) - \rho(t, \mathbf{y}_1) \leq \lVert \mathbf{y}_2 - \mathbf{y}_1 \rVert$.\qed
\end{proof}

The distance function $\rho$ allows the minimal interception time to be expressed as follows:
\begin{equation}
    T^*[\boldsymbol{y}_T] \overset{\mathrm{def}}{=} \min\left\{t \in \mathbb{R}^+_0:\: \rho(t, \boldsymbol{y}_T(t)) \leq \ell\right\} = \inf_{\boldsymbol{u} \in \mathcal{A}} J[\boldsymbol{u}; \boldsymbol{y}_T].
\end{equation}
If $\rho(t, \boldsymbol{y}_T(t)) > \ell$ for all $t \in \mathbb{R}^+_0$, we assume $T^*[\boldsymbol{y}_T] = +\infty$. Thus, the minimal interception time is the first time that the distance between the target position and the reachable set projection does not exceed the capture radius. 

\begin{theorem}\label{th:main_equation}
    For any $\boldsymbol{y}_T \in \mathrm{Lip}_v(\mathbb{R}^+_0, \mathcal{Y})$, if $\rho(0, \boldsymbol{y}_T(0)) \geq \ell$, then
    \begin{equation}
        T^*[\boldsymbol{y}_T] = \min\left\{t \in \mathbb{R}^+_0:\: \rho(t, \boldsymbol{y}_T(t)) = \ell\right\}.
    \end{equation}
\end{theorem}
\begin{proof}
    Using Lemmas~\ref{lem:rho_t_lip},~\ref{lem:rho_y_lip} and the definition of Lipschitz-continuity of $\boldsymbol{y}_T$ for all $t_1, t_2 \in \mathbb{R}^+_0$, we obtain
    \begin{multline*}
        \lvert \rho(t_2, \boldsymbol{y}_T(t_2)) - \rho(t_1, \boldsymbol{y}_T(t_1)) \rvert \leq \lvert t_2 - t_1 \rvert + \lVert \boldsymbol{y}_T(t_2) - \boldsymbol{y}_T(t_1) \rVert \\
        \leq (1 + v) \lvert t_2 - t_1 \rvert.
    \end{multline*}
    It means that $\rho(\cdot, \boldsymbol{y}_T(\cdot))$ is a continuous function of time. Let $\rho(t, \boldsymbol{y}_T(t)) < \ell$ for $t = T^*[\boldsymbol{y}_T] \in \mathbb{R}^+_0$. If $\rho(0, \boldsymbol{y}_T(0)) = \ell$, then $t = 0$. If $\rho(0, \boldsymbol{y}_T(0)) > \ell$, then $t > 0$. According to the intermediate value theorem, there exists $t^* \in (0, t)$ such that $\rho(t^*, \boldsymbol{y}_T(t^*)) = \ell$ and $t^* < t = T^*[\boldsymbol{y}_T]$. That contradiction ends the proof.\qed
\end{proof}

Theorem~\ref{th:main_equation} states that the minimal interception time is the smallest non-negative root of equation $\rho(t, \boldsymbol{y}_T(t)) = \ell$. If $\boldsymbol{y}_T$ is fixed, we can employ the universal root-finding algorithm from Theorem 2 of~\cite{Abaffy2013-rf} to find the least non-negative root of the Lipschitz function $g(\cdot) = \rho(\cdot, \boldsymbol{y}_T(\cdot)) - \ell$. Nevertheless, this algorithm ignores the properties of $\rho$-function, that is, it can be improved. Further, we show how to use these properties to obtain an optimal step size for the corresponding root-finding algorithm.

\subsection{Universal Lower Estimators}

\begin{definition}\label{def:low_est}
    We say that the function $\theta: \mathbb{R}^+_0 \times \mathcal{Y} \to \mathbb{R}^+$ is {\it a universal lower estimator of interception time} if and only if for all $t \in \mathbb{R}^+_0$ and $\mathbf{y} \in \mathcal{Y}$:
    \begin{enumerate}
        \item if $\rho(t, \mathbf{y}) \leq \ell$, then $\theta(t, \mathbf{y}) = t$;
        \item if $\rho(t, \mathbf{y}) > \ell$, then $\theta(t, \mathbf{y}) > t$;
        \item if $\rho(t, \mathbf{y}) > \ell$, then $\rho(\tau, \boldsymbol{y}_T(\tau)) > \ell$ for all $\tau \in [t, \theta(t, \mathbf{y}))$, $\boldsymbol{y}_T \in \mathrm{Lip}_v(\mathbb{R}^+_0, \mathcal{Y})$.
    \end{enumerate}
\end{definition}

This definition is not associated with a distinct form of target trajectory $\boldsymbol{y}_T$. For this reason, the term "universal" is used. Thus, any algorithm for searching for the minimal interception time based on a universal lower estimator is suitable for any target trajectory. The first requirement of the definition is formal, and it sets the value of the lower estimation of the interception time to the given time if point $\mathbf{y}$ is already sufficiently close to the reachable set projection. The second requirement declares a lower estimation to increase if point $\mathbf{y}$ is still sufficiently far from the reachable set projection. According to the third requirement, lower estimations do not overestimate the interception time for any trajectory of target $\boldsymbol{y}_T \in \mathrm{Lip}_v(\mathbb{R}^+_0, \mathcal{Y})$. Later, we provided two universal lower estimators and proved the convergence to the minimal interception time of the corresponding algorithms.

Theoretically, if point $\mathbf{y}$ is greater than $\ell$ away from set $\mathcal{R}(t)$, the third requirement of Definition~\ref{def:low_est} can be satisfied if we consider an estimator that is the infimum of all trajectories of the target for the time taken to approach $\mathcal{R}(t)$ at a distance of no more than $\ell$.
\begin{definition}\label{def:best_est}
    We call {\it the best universal lower estimator} the function given by the rule (Fig.~\ref{fig:best_est}) 
    \begin{equation}
        T(t, \mathbf{y}) \overset{\mathrm{def}}{=} \inf_{\substack{\boldsymbol{y}_T \in \mathrm{Lip}_v(\mathbb{R}^+_0, \mathcal{Y}) \\ \boldsymbol{y}_T(t) = \mathbf{y}}} \min\left\{\tau \in \mathbb{R}^+_0:\: \tau \geq t, \: \rho(\tau, \boldsymbol{y}_T(\tau)) \leq \ell\right\}.
    \end{equation}
    If $\rho(\tau, \boldsymbol{y}_T(\tau)) > \ell$ for all $\tau \in [t, +\infty)$, $\boldsymbol{y}_T \in \mathrm{Lip}_v(\mathbb{R}^+_0, \mathcal{Y})$ such that $\boldsymbol{y}_T(t) = \mathbf{y}$, then we set $T(t, \mathbf{y}) = +\infty$.
\end{definition}

\begin{figure}
    \centering
    \begin{tikzpicture}[scale=1.5]
        \draw[->] (0, -0.25) -- (0, 3) node[above]{$\mathcal{Y}$};
        \draw[->] (0, 0) -- (5, 0);
        \filldraw (0, 0) circle (1pt) node[left]{$0$};
        \draw[dotted] (0, 1.5) -- (3, 1.5);
        \draw[dotted] (3, 0) -- (3, 1.5);
        \filldraw (3, 0) circle (1pt) node[below]{$t$};
        \filldraw (4, 0) circle (1pt) node[below]{$T(t, \mathbf{y})$};
        \filldraw (0, 1.5) circle (1pt) node[left]{$\mathbf{y}$};
        \filldraw (3, 1.5) circle (1pt);
        \draw[dotted] (4, 0) -- (4, 13/12 + 0.22);
        \draw[densely dashdotted] (0, 1) -- (5, 11/6);
        \draw[densely dashdotted] (0, 2) -- (5, 7/6);
        \draw[-] (0, 1.75) .. controls (1, 1.6) and (2, 1.5) .. (3, 1.5) .. controls (3.5, 1.5) and (4.0, 1.4) .. (5, 1.4) node[right]{$\boldsymbol{y}_T(\cdot)$};
        \draw[dashed,yshift=6pt] (0, 0) .. controls (2, 0.7) .. (4, 13/12);
        \draw[-] (0, 0) .. controls (2, 0.7) .. (4, 13/12);
        \draw[-] (0, 0) .. controls (2, 0.5) .. (4, 0.5);
        \draw[dashed,yshift=-6pt] (0, 0) .. controls (2, 0.5) .. (4, 0.5);
        \filldraw (4, 13/12) circle (1pt);
        \filldraw (4, 0.5) circle (1pt) node[above left]{$\mathcal{R}(\cdot)$};
        \filldraw (4, 13/12 + 0.22) circle (1pt);
        \filldraw (4, 0.28) circle (1pt);
        \draw[-] (4, 13/12) -- (4, 0.5);
        \draw[dotted] (0, 13/12) -- (4, 13/12);
        \draw[dotted] (0, 0.5) -- (4, 0.5);
        \filldraw (0, 13/12) circle (1pt);
        \filldraw (0, 0.5) circle (1pt);
        \draw[line width=1] (0, 13/12) -- node[left]{$\mathcal{R}(T(t, \mathbf{y}))$} (0, 0.5);
    \end{tikzpicture}
    \caption{Graphical calculation $T(t, \mathbf{y})$ for given $t$ and $\mathbf{y}$. $\mathcal{Y}$ is depicted as a line for simplicity. $\mathcal{R}$ (solid lines) are surrounded by capture set constraints (dashed lines). Any target trajectory $\boldsymbol{y}_T$ such that $\boldsymbol{y}_T(t) = \mathbf{y}$ lies between the dash-dotted cones because $\boldsymbol{y}_T \in \mathrm{Lip}_v(\mathbb{R}^+, \mathcal{Y})$ (tangent of the cone apex angle is $v$)}
    \label{fig:best_est}
\end{figure}

Clearly, $T$ produces lower estimations of interception time. Moreover, if $\theta$ is an arbitrary universal lower estimator, $\theta(t, \mathbf{y}) \leq T(t, \mathbf{y})$ for any $t \in \mathbb{R}^+_0$ and $\mathbf{y} \in \mathcal{Y}$. Therefore, we can say that $T$ produces lower estimations that are always not less than any other lower estimation of the interception time. Therefore, this function is referred to as the best estimator.

\begin{lemma}\label{lem:T_ineq}
    For all $t \in \mathbb{R}^+_0$ and $\mathbf{y} \in \mathcal{Y}$:
    \begin{equation}
        T(t, \mathbf{y}) = \min\left\{\tau \in \mathbb{R}^+_0:\: \tau \geq t, \: \rho(\tau, \mathbf{y}) \leq v(\tau - t) + \ell\right\}.
    \end{equation}
    If $\rho(\tau, \mathbf{y}) > v(\tau - t) + \ell$ for all $\tau \in [t, +\infty)$, then we assume that the right part of the equation is equal to $+\infty$.
\end{lemma}
\begin{proof}
    Let us first prove that $T(t, \mathbf{y})$ is equal to
    \begin{equation}\label{eq:T_tmp_form}
        \min_{\boldsymbol{\eta} \in \mathcal{Y}}\min\left\{\tau \in \mathbb{R}^+_0:\: \tau \geq t, \: \rho(\tau, \boldsymbol{\eta}) \leq \ell, \: \lVert \boldsymbol{\eta} - \mathbf{y} \rVert \leq v(\tau - t) \right\}.
    \end{equation}
    Let $\boldsymbol{\eta} \in \mathcal{Y}$. Set 
    \begin{equation}
        \boldsymbol{y}_T(\tau) = 
        \begin{cases}
            \mathbf{y} + v\frac{\boldsymbol{\eta} - \mathbf{y}}{\lVert \boldsymbol{\eta} - \mathbf{y} \rVert}(\tau - t), \quad \boldsymbol{\eta} \neq \mathbf{y};\\
            \mathbf{y}, \quad \boldsymbol{\eta} = \mathbf{y}.
        \end{cases}
    \end{equation}
    Since $\boldsymbol{y}_T \in \mathrm{Lip}_v(\mathbb{R}^+_0, \mathcal{Y})$ and $\boldsymbol{y}_T(t) = \mathbf{y}$, the value $T(t, \mathbf{y})$ given by Definition~\ref{def:best_est} is less than or equal to~\eqref{eq:T_tmp_form}. On the other hand, $\rho(\tau, \boldsymbol{y}_T(\tau)) \leq \ell$ follows from $\rho(\tau, \boldsymbol{\eta}) \leq \ell$ and $\lVert \boldsymbol{\eta} - \mathbf{y} \rVert \leq v(\tau - t)$ if we fix $\boldsymbol{\eta} = \boldsymbol{y}_T(\tau)$. Thus, $T(t, \mathbf{y})$ given by Definition~\ref{def:best_est} is equal to or greater than~\eqref{eq:T_tmp_form}. Next, for any $\boldsymbol{\eta} \in \mathcal{Y}$ such that $\lVert \boldsymbol{\eta} - \mathbf{y} \rVert \leq v(\tau - t)$ and $\rho(\tau, \boldsymbol{\eta}) \leq \ell$,
    \begin{equation}
        \rho(\tau, \mathbf{y}) \leq \min_{\tilde{\boldsymbol{\eta}} \in \mathcal{R}(\tau)} \lVert \tilde{\boldsymbol{\eta}} - \boldsymbol{\eta} \rVert + \lVert \boldsymbol{\eta} - \mathbf{y} \rVert = \rho(\tau, \boldsymbol{\eta}) + \lVert \boldsymbol{\eta} - \mathbf{y} \rVert \leq v(\tau - t) + \ell.
    \end{equation}
    By the above,
    \begin{equation}
        T(t, \mathbf{y}) \geq \min\left\{\tau \in \mathbb{R}^+_0:\: \tau \geq t, \: \rho(\tau, \mathbf{y}) \leq v(\tau - t) + \ell\right\}.
    \end{equation}
    We fix $\tilde{\boldsymbol{\eta}} \in \mathrm{arg}\min_{\boldsymbol{\eta} \in \mathcal{R}(\tau)} \lVert \boldsymbol{\eta} - \mathbf{y} \rVert$. If $\lVert \tilde{\boldsymbol{\eta}} - \mathbf{y} \rVert < \ell$, we set $\boldsymbol{\eta} = \mathbf{y}$. It is easy to verify that $\rho(\tau, \boldsymbol{\eta}) \leq \ell$ and $\lVert \boldsymbol{\eta} - \mathbf{y} \rVert \leq v(\tau - t)$ if $\tau \geq t$. If $\lVert \tilde{\boldsymbol{\eta}} - \mathbf{y} \rVert \geq \ell$, then we consider 
    \begin{equation}
        \boldsymbol{\eta} = \mathbf{y} + (\tilde{\boldsymbol{\eta}} - \mathbf{y})\left(1 - \frac{\ell}{\lVert \tilde{\boldsymbol{\eta}} - \mathbf{y} \rVert}\right).
    \end{equation}
    Now, if $\tau \geq t$ and $\rho(\tau, \mathbf{y}) \leq v(\tau - t) + \ell$, then
    \begin{equation}
        \lVert \boldsymbol{\eta} - \mathbf{y} \rVert = \lVert \tilde{\boldsymbol{\eta}} - \mathbf{y} \rVert - \ell = \rho(\tau, \mathbf{y}) - \ell \leq v(\tau - t),
    \end{equation}
    and
    \begin{equation}
        \rho(\tau, \boldsymbol{\eta}) = \min_{\tilde{\mathbf{y}} \in \mathcal{R}(\tau)} \lVert \boldsymbol{\eta} - \tilde{\mathbf{y}} \rVert \leq \lVert \boldsymbol{\eta} - \tilde{\boldsymbol{\eta}} \rVert = \ell.
    \end{equation}
    Thus,
    \begin{equation}
        T(t, \mathbf{y}) \leq \min\left\{\tau \in \mathbb{R}^+_0:\: \tau \geq t, \: \rho(\tau, \mathbf{y}) \leq v(\tau - t) + \ell\right\}.
    \end{equation}\qed
\end{proof}

We now clarify how to calculate $T(t, \mathbf{y})$ as the root of a real equation.
\begin{theorem}\label{th:T_calc}
    If $t \in \mathbb{R}^+_0$, $\mathbf{y} \in \mathcal{Y}$ such that $\rho(t, \mathbf{y}) \geq \ell$, then
    \begin{equation}
        T(t, \mathbf{y}) = \min\left\{\tau \in \mathbb{R}^+_0:\: \tau \geq t, \: \rho(\tau, \mathbf{y}) = v(\tau - t) + \ell\right\}.
    \end{equation}
    If $\rho(\tau, \mathbf{y}) > v(\tau - t) + \ell$ for all $\tau \in [t, +\infty)$, then we suppose that the right part of the equation is equal to $+\infty$.
\end{theorem}
\begin{proof}
    Case $\rho(t, \mathbf{y}) = \ell$ is trivial. Suppose $\rho(t, \mathbf{y}) > \ell$. Let $\tau \in [t, T(t, \mathbf{y}))$. Using Lemma~\ref{lem:T_ineq}, we deduce $\rho(\tau, \mathbf{y}) > v(\tau - t) + \ell$. Since $\rho(\cdot, \mathbf{y})$ is a continuous function, for $\tau \to T(t, \mathbf{y}) - 0$ we have $\rho(T(t, \mathbf{y}), \mathbf{y}) \geq v(T(t, \mathbf{y}) - t) + \ell$. Lemma~\ref{lem:T_ineq} also gives $\rho(T(t, \mathbf{y}), \mathbf{y}) \leq v(T(t, \mathbf{y}) - t) + \ell$.\qed
\end{proof}

$T$ is the best universal lower estimator, but the calculation of $T(t, \mathbf{y})$ requires knowledge of the solution of $\rho(\tau, \mathbf{y}) = v(\tau - t) + \ell$ which is not always available in practice. Further, we describe a universal lower estimator that can be easily calculated if $\rho(t, \mathbf{y})$ is known.

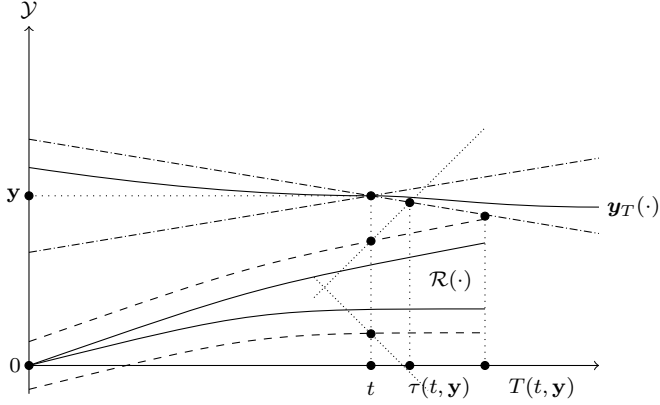
\begin{figure}
    \centering
    \begin{tikzpicture}[scale=1.5]
        \draw[->] (0, -0.25) -- (0, 3) node[above]{$\mathcal{Y}$};
        \draw[->] (0, 0) -- (5, 0);
        \filldraw (0, 0) circle (1pt) node[left]{$0$};
        \draw[dotted] (0, 1.5) -- (3, 1.5);
        \draw[dotted] (3, 0) -- (3, 1.5);
        \filldraw (3, 0) circle (1pt);
        \filldraw (4, 0) circle (1pt);
        \filldraw (0, 1.5) circle (1pt) node[left]{$\mathbf{y}$};
        \filldraw (3, 1.5) circle (1pt);
        \draw[densely dashdotted] (0, 1) -- (5, 11/6);
        \draw[densely dashdotted] (0, 2) -- (5, 7/6);
        \draw[-] (0, 1.75) .. controls (1, 1.6) and (2, 1.5) .. (3, 1.5) .. controls (3.5, 1.5) and (4.0, 1.4) .. (5, 1.4) node[right]{$\boldsymbol{y}_T(\cdot)$};
        \node at (3.7, 0.75) {$\mathcal{R}(\cdot)$};
        \draw[dashed,yshift=6pt] (0, 0) .. controls (2, 0.7) .. (4, 13/12);
        \draw[-] (0, 0) .. controls (2, 0.7) .. (4, 13/12);
        \draw[-] (0, 0) .. controls (2, 0.5) .. (4, 0.5);
        \draw[dashed,yshift=-6pt] (0, 0) .. controls (2, 0.5) .. (4, 0.5);
        \filldraw (3, 1.1) circle (1pt);
        \filldraw (3.34, 1.44) circle (1pt);
        \filldraw (4, 1.32) circle (1pt);
        \draw[densely dotted] (2.5, 0.6) -- (4, 2.1);
        \draw[densely dotted] (2.5, 0.78) -- (3.5, -0.22);
        \filldraw (3, 0.28) circle (1pt);
        \draw[dotted] (3.34, 1.44) -- (3.34, 0);
        \draw[dotted] (4, 0) -- (4, 1.32);
        \filldraw (3.34, 0) circle (1pt);
        \node at (3.0, -0.2) {$t$};
        \node at (3.6, -0.2) {$\tau(t, \mathbf{y})$};
        \node at (4.5, -0.2) {$T(t, \mathbf{y})$};
    \end{tikzpicture}
    \caption{Graphical calculation $\tau(t, \mathbf{y})$ for a given $t$ and $\mathbf{y}$. $\mathcal{Y}$ is depicted as a line for simplicity. $\mathcal{R}$ (solid lines) are surrounded by capture set constraints (dashed lines). Any target trajectory $\boldsymbol{y}_T$ such that $\boldsymbol{y}_T(t) = \mathbf{y}$ lies between the dash-dotted lines, because $\boldsymbol{y}_T \in \mathrm{Lip}_v(\mathbb{R}^+, \mathcal{Y})$ (the tangent of the cone apex angle is $v$). The densely dotted cone has the tangent of the apex angle, which is equal to $1$}
    \label{fig:simple_est}
\end{figure}

\begin{definition}\label{def:tau}
    We call {\it the simple universal lower estimator} the function given by the rule (Fig.~\ref{fig:simple_est})
    \begin{equation}
        \tau(t, \mathbf{y}) \overset{\mathrm{def}}{=}
        \begin{cases}
            t + \frac{\rho(t, \mathbf{y}) - \ell}{1 + v}, \quad \rho(t, \mathbf{y}) > \ell;\\
            t, \quad \rho(t, \mathbf{y}) \leq \ell.
        \end{cases}
    \end{equation}
\end{definition}

The proposed universal lower estimator uses the geometrical restriction: $\lvert \boldsymbol{f}(\mathbf{x}, \mathbf{u}) \rvert \leq 1$. This design of the function implies the movement of the plant and target towards each other at speeds of $1$ and $v$, respectively. Further, we prove that this function satisfies the requirements of Definition~\ref{def:low_est}.

\begin{lemma}\label{lem:tau_less_eq_T}
    $\tau(t, \mathbf{y}) \leq T(t, \mathbf{y})$ for all $t \in \mathbb{R}^+_0$, $\mathbf{y} \in \mathcal{Y}$.
\end{lemma}
\begin{proof}
    Using $T(t, \mathbf{y}) \geq t$ and Lemma~\ref{lem:rho_t_lip}, we obtain
    \begin{equation}\label{eq:tmp_lem41}
        \rho(t, \mathbf{y}) - \rho(T(t, \mathbf{y}), \mathbf{y}) \leq T(t, \mathbf{y}) - t.
    \end{equation}
    According to Theorem~\ref{th:T_calc},
    \begin{equation}\label{eq:tmp_lem42}
        \rho(T(t, \mathbf{y}), \mathbf{y}) = v(T(t, \mathbf{y}) - t) + \ell.
    \end{equation}
    Substituting~\eqref{eq:tmp_lem42} and $\rho(t, \mathbf{y}) = (1 + v)(\tau(t, \mathbf{y}) - t) + \ell$ to~\eqref{eq:tmp_lem41}
    we obtain $\tau(t, \mathbf{y}) \leq T(t, \mathbf{y})$.\qed
\end{proof}

\begin{lemma}\label{lem:tau_is_low_est}
    $\tau$ is a universal lower estimator of the interception time.
\end{lemma}
\begin{proof}
    1. If $\rho(t, \mathbf{y}) \leq \ell$, then $\tau(t, \mathbf{y}) = t$. 2. If $\rho(t, \mathbf{y}) > \ell$, then
    \begin{equation}
        \tau(t, \mathbf{y}) - t = \frac{\rho(t, \mathbf{y}) - \ell}{1 + v} > \frac{\ell - \ell}{1 + v} = 0.
    \end{equation}
    3. Since $T$ is a universal lower estimator of the interception time, Lemma~\ref{lem:tau_less_eq_T} guarantees that the corresponding property holds.\qed
\end{proof}

\subsection{Iterative Algorithms}
Here, we define fixed-point iteration algorithms that use universal lower estimators of interception time.
\begin{definition}
    An iterative algorithm based on a universal lower estimator $\theta$ is {\it correctly} convergent if the sequence $\{t_n\}_{n=0}^\infty$ such that $t_0 = 0$ and $t_n = \theta(t_{n - 1}, \boldsymbol{y}_T(t_{n - 1}))$ converges to $T^*[\boldsymbol{y}_T]$ for any $\boldsymbol{y}_T \in \mathrm{Lip}_v(\mathbb{R}^+_0, \mathcal{Y})$.
\end{definition}

\begin{lemma}\label{lem:tau_is_conv}
    The iterative algorithm based on $\tau$ is correctly convergent.
\end{lemma}
\begin{proof}
    Fix $\boldsymbol{y}_T \in \mathrm{Lip}_v(\mathbb{R}^+_0, \mathcal{Y})$. According to Definition~\ref{def:tau},
    \begin{equation}
        t_n - t_{n - 1} =
        \begin{cases}
            \frac{\rho(t_{n - 1}, \boldsymbol{y}_T(t_{n - 1})) - \ell}{1 + v}, \quad \rho(t_{n - 1}, \boldsymbol{y}_T(t_{n - 1})) > \ell;\\
            0, \quad \rho(t_{n - 1}, \boldsymbol{y}_T(t_{n - 1})) \leq \ell.
        \end{cases}
    \end{equation}
    Thus, $\{t_n\}_{n = 0}^\infty$ is a non-decreasing sequence. We first consider the case $T^*[\boldsymbol{y}_T] < +\infty$. Suppose that there exists $n \in \mathbb{N}$ such that $t_n > T^*[\boldsymbol{y}_T]$. Without loss of generality we assume that $n$ is minimal one and $t_{n - 1} \leq T^*[\boldsymbol{y}_T]$. If $t_{n - 1} = T^*[\boldsymbol{y}_T]$, then $\rho(t_{n - 1}, \boldsymbol{y}_T(t_{n - 1})) \leq \ell$ and $t_n = \rho(t_{n - 1}, \boldsymbol{y}_T(t_{n - 1})) = t_{n - 1} = T^*[\boldsymbol{y}_T]$. If $t_{n - 1} < T^*[\boldsymbol{y}_T]$, then $\rho(t_{n - 1}, \boldsymbol{y}_T(t_{n - 1})) > \ell$. Lemma~\ref{lem:tau_is_low_est} states that $\tau$ is a lower estimation of interception time. According to requirement 3 of Definition~\ref{def:low_est}, $t_n = \tau(t_{n - 1}, \boldsymbol{y}_T(t_{n - 1}))$ cannot be greater than $T^*[\boldsymbol{y}_T]$. Thus, $t_n \leq T^*[\boldsymbol{y}_T]$ and
    \begin{equation}
        \lim_{n \to \infty} t_n = t^* \leq T^*[\boldsymbol{y}_T].
    \end{equation}
    Let $t^* < T^*[\boldsymbol{y}_T]$. Since $T^*[\boldsymbol{y}_T]$ is the minimal interception time, $\rho(t^*, \boldsymbol{y}_T(t^*)) - \ell > 0$. Let $\mu = \rho(t^*, \boldsymbol{y}_T(t^*)) - \ell$. It follows from the continuity of $\rho(\cdot, \boldsymbol{y}_T(\cdot))$ that there exists $N \in \mathbb{N}$ such that for all $n > N$:
    \begin{equation}
        \rho(t^*, \boldsymbol{y}_T(t^*)) - \rho(t_{n - 1}, \boldsymbol{y}_T(t_{n - 1})) < \frac\mu2.
    \end{equation}
    It implies
    \begin{equation}
        \mu = \rho(t^*, \boldsymbol{y}_T(t^*)) - \ell < \rho(t_{n - 1}, \boldsymbol{y}_T(t_{n - 1})) + \frac\mu2 - \ell = (1 + v)(t_n - t_{n - 1}) + \frac\mu2
    \end{equation}
    and $2(t_n - t_{n - 1}) > \mu/(1 + v)$ for all $n > N$. It contradicts the convergence of $\{t_n\}_{n = 0}^\infty$. 
    
    Now we turn to case $T^*[\boldsymbol{y}_T] = +\infty$. Further we will show that for any $t^* \in \mathbb{R}^+_0$ there exists $n \in \mathbb{N}$ such that $t_n > t^*$, that is,
    \begin{equation}
        \lim_{n \to \infty} t_n = T^*[\boldsymbol{y}_T] = +\infty.
    \end{equation}
    Fix $t^* \in \mathbb{R}^+_0$ and define 
    \begin{equation}
        \mu \overset{\mathrm{def}}{=} \min_{t \in [0, t^*]} \rho(t, \boldsymbol{y}_T(t)) - \ell.
    \end{equation}
    We conclude from $T^*[\boldsymbol{y}_T] \notin [0, t^*]$ that $\mu > 0$, hence that
    \begin{equation}
        t_n - t_{n - 1} = \frac{\rho(t_{n - 1}, \boldsymbol{y}_T(t_{n - 1})) - \ell}{1 + v} \geq \frac{\mu}{1 + v} > 0,
    \end{equation}
    and finally that
    \begin{equation}
        t_n \geq \frac{n\mu}{1 + v}.
    \end{equation}\qed
\end{proof}

\begin{theorem}\label{th:T_corr_conv}
    The iterative algorithm based on $T$ is correctly convergent.
\end{theorem}
\begin{proof}
    It follows from Lemmas~\ref{lem:tau_less_eq_T},~\ref{lem:tau_is_conv} and Definition~\ref{def:best_est}.\qed
\end{proof}

Thus, $\tau$ and $T$ can be used to compute $T^*[\boldsymbol{y}_T]$ using a fixed-point iteration algorithm. If $T(t, \mathbf{y})$ can be obtained from Theorem~\ref{th:T_calc}, then it is appropriate to use $T(t, \mathbf{y})$ instead $\tau(t, \mathbf{y})$ because the corresponding step size $t_n - t_{n - 1}$ is maximal among the convergent algorithms that use only Lipschitz properties of $\boldsymbol{y}_T$.

The simple universal lower estimator helps solve the minimal root-finding problem for $\rho(t, \boldsymbol{y}_T(t)) = \ell$ in the manner described in~\cite{Abaffy2013-rf}. The novelty of the current study lies in the use of additional information about the plant (the distance function $\rho$) using the best universal lower estimator $T$ to increase the guaranteed step of the fixed-point iteration algorithm.

The rate of convergence of the fixed-point iteration algorithm using function $T$ can be arbitrarily small or high for specific cases of trajectory $\boldsymbol{y}_T$. If the target moves to the nearest point on the reachable set in a straight line, then the algorithm converges in one step. To demonstrate that the convergence rate of the algorithm can be arbitrarily small, we consider a specific example. Let $\mathcal{X} = \mathcal{Y} = \mathbb{R}$. The trajectory of the target is given by $\boldsymbol{y}_T(t) = \ell + v\alpha(1 - t)$, where $\alpha \in \mathrm{Lip}_1(\mathbb{R}^+_0, \mathbb{R}^+_0)$ and $\alpha(0) = 0$, $\alpha(1 - t) > 0$ for all $t \in [0, 1)$. Let $\boldsymbol{f}(\mathbf{x}, \mathbf{u}) = 0$ and $\mathbf{x}_0 = 0$. We can derive that $T(t, \mathbf{y}) = t + (\mathbf{y} - \ell)/v$. Thus, $t_n = t_{n - 1} + \alpha(1 - t_{n - 1})$. Step size $t_n - t_{n - 1} = \alpha(1 - t_{n - 1})$ can be arbitrarily small because there are no additional restrictions on $\alpha$. The same conclusions can be drawn regarding any other correctly convergent algorithm based on a universal lower estimator of interception time. Thus, the cost of universality of the proposed method is ineffectiveness for some cases of target motion. However, in practical terms, these cases are nothing more than mathematical anomalies, because the target must smoothly touch the boundary of the reachable set in these cases.

\begin{figure}
    \centering
    \begin{tikzpicture}[scale=1.5]
        \draw[->] (0, -0.45) -- (0, 3) node[above]{$\mathcal{Y}$};
        \draw[->] (0, 0) -- (5, 0);
        \filldraw (0, 0) circle (1pt) node[left]{$0$};
        \draw[-] (0, 1.75) .. controls (1, 1.6) and (2, 1.5) .. (3, 1.5) .. controls (3.5, 1.5) and (4.0, 1.4) .. (5, 1.4);
        \node at (1.5, 1.8) {$\boldsymbol{y}_T(\cdot)$};
        \node at (3.2, 0.70) {$\mathcal{R}(\cdot)$};
        \draw[densely dotted,yshift=10pt] (0, 0) .. controls (2, 0.7) .. (4, 13/12) .. controls (4.5, 1.2) .. (5, 1.3);
        \draw[dashed,yshift=6pt] (0, 0) .. controls (2, 0.7) .. (4, 13/12) .. controls (4.5, 1.2) .. (5, 1.3);
        \draw[-] (0, 0) .. controls (2, 0.7) .. (4, 13/12) .. controls (4.5, 1.2) .. (5, 1.3);
        \draw[-] (0, 0) .. controls (2, 0.5) .. (4, 0.5) .. controls (4.5, 0.47) .. (5, 0.42);
        \draw[dashed,yshift=-6pt] (0, 0) .. controls (2, 0.5) .. (4, 0.5) .. controls (4.5, 0.47) .. (5, 0.42);
        \draw[densely dotted,yshift=-10pt] (0, 0) .. controls (2, 0.5) .. (4, 0.5) .. controls (4.5, 0.47) .. (5, 0.42);
        \filldraw (0, 1.75) circle (1pt) node[left]{$\boldsymbol{y}_T(0)$};
        \draw[dashdotted] (0, 1.75) -- (3, 1);
        \draw[dashdotted] (0, 1.75) -- (3, 2.5);
        \filldraw (2.78, 1.05) circle (1pt);
        \filldraw (2.78, 0) circle (1pt) node[below]{$t_1$};
        \draw[dotted] (2.78, 1.5) -- (2.78, 0);
        \filldraw (2.78, 1.5) circle (1pt);
        \node at (2.78, 1.7) {$\boldsymbol{y}_T(t_1)$};
        \draw[dashdotted] (2.78, 1.5) -- (3.78, 1.25);
        \draw[dashdotted] (2.78, 1.5) -- (3.78, 1.75);
        \filldraw (3.76, 1.25) circle (1pt);
        \filldraw (3.76, 0) circle (1pt) node[below]{$t_2$};
        \filldraw (3.76, 1.45) circle (1pt);
        \node at (3.76, 1.9) {$\boldsymbol{y}_T(t_2)$};
        \draw[dotted] (3.76, 1.45) -- (3.76, 0);
        \draw[dashdotted] (3.76, 1.45) -- (4.76, 1.25);
        \draw[dashdotted] (3.76, 1.45) -- (4.76, 1.70);
        \filldraw (4.25, 1.35) circle (1pt);
        \filldraw (4.25, 0) circle (1pt) node[below]{$t_3$};
        \filldraw (4.25, 1.42) circle (1pt);
        \node at (4.45, 1.8) {$\boldsymbol{y}_T(t_3)$};
        \draw[dotted] (4.25, 0) -- (4.25, 1.42);
        \filldraw (4.48, 1.41) circle (1pt);
        \filldraw (4.48, 0) circle (1pt);
        \node at (4.8, -0.15) {$T^*[\boldsymbol{y}_T]$};
        \draw[dotted] (4.48, 0) -- (4.48, 1.41);
    \end{tikzpicture}
    \caption{Iterations of Algorithm~\ref{alg:comp} in the case $\theta = T$ for a given $\boldsymbol{y}_T$. $\mathcal{Y}$ is depicted as a line for simplicity. $\mathcal{R}$ (solid lines) are surrounded by capture set constraints (dashed and densely dotted lines). The dashed lines correspond to the desired capture radius $\ell$, while the densely dotted lines denote the capture radius $\ell(1 + \varepsilon)$. The dash-dotted cones have the tangent of the cone apex angle that is equal to $v$. The iterations stopped in the third step because $\rho(t_3, \boldsymbol{y}_T(t_3)) \leq \ell(1 + \varepsilon)$ (the intersection point appears below the densely dotted line)}
    \label{fig:one_step}
\end{figure}
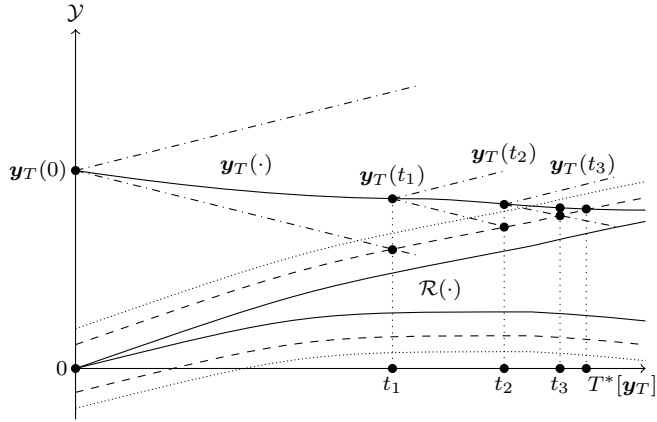

Based on these results, we propose a simple algorithm for computing an approximation of the minimal interception time (Algorithm~\ref{alg:comp}). The given value of $\varepsilon$ sets the admissible relative error when approaching the capture set of a moving target. We used this value to limit the number of steps to a finite number. If the value of $T$ can be calculated without difficulty, then it is better to pass $T$ instead of $\theta$. In other cases, $\tau$ is used instead of $\theta$. An example of this algorithm work in the case $\theta = T$ is shown in Fig.~\ref{fig:one_step}.

\begin{algorithm}
    \caption{Computing an approximation of $T^*[\boldsymbol{y}_T]$}\label{alg:comp}
    \textbf{Given:} $\ell$ (capture radius), $\varepsilon$ (relative error), $v$ (the maximal target speed), $\rho(\cdot, \cdot)$ (distance to the reachable set projection), $\theta(\cdot, \cdot)$ (universal lower estimator $T(\cdot, \cdot)$ or $\tau(\cdot, \cdot)$)\\
    \textbf{Input: } $\boldsymbol{y}_T(\cdot)$ (trajectory of the moving target)
    \begin{algorithmic}
        \Require $\ell > 0$, $\varepsilon > 0$, $v \geq 0$, $\boldsymbol{y}_T \in \mathrm{Lip}_v(\mathbb{R}, \mathcal{Y})$
        \State $t \gets 0$
        \While{$\rho(t, \boldsymbol{y}_T(t)) > \ell(1 + \varepsilon)$}
            \State $t \gets \theta(t, \boldsymbol{y}_T(t))$
        \EndWhile
    \end{algorithmic}
    \textbf{Output: } $t$ (approximation of the minimal interception time)
\end{algorithm}

\section{Examples}\label{sec:examples}
In this section, we provide two examples to demonstrate the proposed algorithmic framework. In both cases, the path of the target lies on plane $\mathcal{Y} = \mathbb{R}^2$. Throughout the section, $\mathbf{y} = \begin{bmatrix} x & y \end{bmatrix}^\top$ denotes a point on $\mathcal{Y}$. The norm of $\mathcal{Y}$ is given by $\lVert \mathbf{y} \rVert = \sqrt{x^2 + y^2}$.

\subsection{Simple Motions}
The state space $\mathcal{X}$ coincides with the plane $\mathcal{Y}$ and $\mathcal{Z} = \varnothing$. The set of admissible values of the control is $\mathcal{U} = \{\mathbf{u} \in \mathcal{Y}: \lVert \mathbf{u} \rVert \leq 1\}$. The state equation is the following: $\boldsymbol{f}(\mathbf{x}, \mathbf{u}) = \mathbf{u}$.

The reachable set of the system at time $t \in \mathbb{R}^+_0$ is a disk with the center at $\boldsymbol{0} \in \mathcal{Y}$. The radius of the disk is $t$. The distance from $\mathbf{y} \in \mathcal{Y}$ to the disk is given by
\begin{equation}
    \rho(t, \mathbf{y}) =
    \begin{cases}
        0, \quad \lVert \mathbf{y} \rVert < t;\\
        \lVert \mathbf{y} \rVert - t, \quad \lVert \mathbf{y} \rVert \geq t.\\
    \end{cases}
\end{equation}
Thus, Definition~\ref{def:tau} can be rewritten as
\begin{equation}
    \tau(t, \mathbf{y}) =
    \begin{cases}
        \frac{\lVert \mathbf{y} \rVert + vt - \ell}{1 + v}, \quad \lVert \mathbf{y} \rVert > t + \ell;\\
        t, \quad \lVert \mathbf{y} \rVert \leq t + \ell.
    \end{cases}
\end{equation}
The real equation of Theorem~\ref{th:T_calc} can be explicitly solved. We can easily show that $T(t, \mathbf{y}) = \tau(t, \mathbf{y})$. According to Theorem~\ref{th:T_corr_conv}, sequence $\{t_n\}_{n = 0}^\infty$ given by $t_0 = 0$ and $t_n = \tau(t_{n - 1}, \boldsymbol{y}_T(t_{n - 1}))$ converges to $T^*[\boldsymbol{y}_T]$ for any $\boldsymbol{y}_T \in \mathrm{Lip}_v(\mathbb{R}^+_0, \mathcal{Y})$. Moreover, the step size $t_n - t_{n - 1}$ is maximal, because $\tau = T$. In Section~\ref{sec:numerical_experiments}, we provide numerical experiments for this case.

\subsection{Dubins Car}
The state space for the Dubins model is $\mathcal{X} = \mathbb{R}^2 \times \mathbb{S}$. Here, $\mathcal{Z} = \mathbb{S}$ is the space of angles. The state vector is $\boldsymbol{x}(t) = \begin{bmatrix} x(t) & y(t) & \varphi(t) \end{bmatrix}^\top$. The control input is given by the scalar $\boldsymbol{u}(t) = u(t) \in \mathcal{U} = [-1, +1]$. The state equation is the following:
\begin{equation}
    \boldsymbol{f}(\boldsymbol{x}, \boldsymbol{u}) =
    \begin{bmatrix}
        \cos\varphi\\
        \sin\varphi
    \end{bmatrix}, \quad \boldsymbol{g}(\boldsymbol{x}, \boldsymbol{u}) = u.
\end{equation}
Without loss of generality, we set $\mathbf{x}_0 = \begin{bmatrix} 0 & 0 & \pi/2 \end{bmatrix}^\top$. In this section, we use the common notation and some known results related to the Dubins model. Similar to~\cite{Boissonnat1994-pd} we assume that $C$ denotes an arc of a unit circle and $S$ denotes a straight-line segment. Sequences $CS$ and $CC$ are assigned to a trajectory consisting of only $C$ and $S$ consecutive parts. If $C$ corresponds to a clockwise or counterclockwise turn, it is replaced by $R$ or $L$ respectively. For the following purposes, we used a slightly modified partition of $\mathcal{Y}$ from~\cite{Ding2019-uw}:
\begin{align*}
    &\mathcal{D}_I \overset{\mathrm{def}}{=} \left\{\mathbf{y} \in \mathcal{Y}:\: \alpha_{CS}(\mathbf{y}) < 0 \vee \mathbf{y} = \boldsymbol{0}\right\},\\
    &\mathcal{D}_{II} \overset{\mathrm{def}}{=} \left\{\mathbf{y} \in \mathcal{Y}:\: \mathbf{y} \notin \mathcal{D}_I, \: \mathbf{y} \notin \mathcal{D}_{III}\right\},\\
    &\mathcal{D}_{III} \overset{\mathrm{def}}{=} \left\{\mathbf{y} \in \mathcal{Y}:\: \alpha_{CC}(\mathbf{y}) > -1, \: y > 0, \: \mathbf{y} \notin \mathcal{D}_I\right\}.
\end{align*}
Here,
\begin{equation}
    \alpha_{CS}(\mathbf{y}) \overset{\mathrm{def}}{=} (1 - \lvert x \rvert)^2 + y^2 - 1, \quad \alpha_{CC}(\mathbf{y}) \overset{\mathrm{def}}{=} \frac{5 - (1 + \lvert x \rvert)^2 - y^2}4.
\end{equation}
Eqs. (4.3), (4.6) from~\cite{Buzikov2021-md} define lengths of Dubins' paths that lead to the boundary of planar reachable set. We will denote by $V_{CS}(\mathbf{y})$, $V_{CC}^+(\mathbf{y})$, $V_{CC}^-(\mathbf{y})$ these lengths. Gathering results from~\cite{Buzikov2021-md}, we set
\begin{align*}
    &V_{CS}(\mathbf{y}) \overset{\mathrm{def}}{=} \theta_{CS}(\mathbf{y}) + \sqrt{\alpha_{CS}(\mathbf{y})}, \quad \mathbf{y} \in \mathcal{D}_{II} \cup \mathcal{D}_{III} \text{ or } (\mathbf{y} = \boldsymbol{0});\\
    &V_{CC}^+(\mathbf{y}) \overset{\mathrm{def}}{=} \theta_{CC}^+(\mathbf{y}) + \arccos\alpha_{CC}(\mathbf{y}), \quad \mathbf{y} \in \mathcal{D}_{III};\\
    &V_{CC}^-(\mathbf{y}) \overset{\mathrm{def}}{=} \theta_{CC}^-(\mathbf{y}) + 2\pi - \arccos\alpha_{CC}(\mathbf{y}), \quad \mathbf{y} \in \mathcal{D}_I \cup \mathcal{D}_{III}.
\end{align*}
Here,
\begin{align*}
    &\theta_{CS}(\mathbf{y}) =
    \begin{cases}
        \arccos\frac{1 - \lvert x \rvert + y\sqrt{\alpha_{CS}(\mathbf{y})}}{1 + \alpha_{CS}(\mathbf{y})}, \quad y \geq (1 - \lvert x \rvert)\sqrt{\alpha_{CS}(\mathbf{y})};\\
        2\pi - \arccos\frac{1 - \lvert x \rvert + y\sqrt{\alpha_{CS}(\mathbf{y})}}{1 + \alpha_{CS}(\mathbf{y})}, \quad y < (1 - \lvert x \rvert)\sqrt{\alpha_{CS}(\mathbf{y})},
    \end{cases}\\
    &\theta_{CC}^\pm(\mathbf{y}) = \arccos\frac{(1 + \lvert x \rvert)(2 - \alpha_{CC}(\mathbf{y})) \pm y\sqrt{1 - \alpha_{CC}^2(\mathbf{y})}}{(1 + \lvert x \rvert)^2 + y^2}.
\end{align*}

The planar reachable set of the Dubins car~\cite{Cockayne1975-pu,Buzikov2021-md} can be described in the explicit form:
\begin{multline*}
    \mathcal{R}(t) = \left\{\mathbf{y} \in \mathcal{Y}: \left(\mathbf{y} \in \mathcal{D}_{III}, \: t \geq V_{CS}(\mathbf{y}), \: (t \geq V_{CC}^-(\mathbf{y}) \vee V_{CC}^+(\mathbf{y}) \geq t)\right)\right.\\
    \left. \vee \left(\mathbf{y} \in \mathcal{D}_{II}, \: t \geq V_{CS}(\mathbf{y})\right) \vee \left(\mathbf{y} \in \mathcal{D}_I, \: t \geq V_{CC}^-(\mathbf{y})\right) \vee \left(t = 0, \: \mathbf{y} = \boldsymbol{0}\right) \right\}.
\end{multline*}
$\mathcal{R}(t)$ is a closed set for any $t \in \mathbb{R}^+_0$. Hence, if $\mathbf{y} \in \mathcal{Y} \setminus \mathcal{R}(t)$, then the nearest point to $\mathbf{y}$ lies on the boundary of $\mathcal{R}(t)$. As we know from~\cite{Cockayne1975-pu,Buzikov2021-md}, the boundary of the planar reachable set consists of two parts $\mathcal{B}_{CS}(t)$, $\mathcal{B}_{CC}(t)$. Thus, the distance from $\mathbf{y} \in \mathcal{Y} \setminus \mathcal{R}(t)$ to $\mathcal{R}(t)$ can be calculated in the following way:
\begin{equation}
    \rho(t, \mathbf{y}) = \min_{\boldsymbol{\eta} \in \mathcal{R}(t)} \lVert \boldsymbol{\eta} - \mathbf{y} \rVert = \min_{\boldsymbol{\eta} \in \partial\mathcal{R}(t)} \lVert \boldsymbol{\eta} - \mathbf{y} \rVert
    = \min(\rho_{CS}(t, \mathbf{y}), \rho_{CC}(t, \mathbf{y})),
\end{equation}
where
\begin{equation}
    \rho_{CS}(t, \mathbf{y}) \overset{\mathrm{def}}{=} \min_{\boldsymbol{\eta} \in \mathcal{B}_{CS}(t)} \lVert \boldsymbol{\eta} - \mathbf{y} \rVert, \quad \rho_{CC}(t, \mathbf{y}) \overset{\mathrm{def}}{=} \min_{\boldsymbol{\eta} \in \mathcal{B}_{CC}(t)} \lVert \boldsymbol{\eta} - \mathbf{y} \rVert.
\end{equation}
\begin{assumption}\label{assum:rho}
    For any $t \in \mathbb{R}^+_0$, if $\mathbf{y} \in \mathcal{Y} \setminus \mathcal{R}(t)$, then
    \begin{equation}
        \rho(t, \mathbf{y}) =
        \begin{cases}
            V_{CS}(\mathbf{y}) - t, \quad (\mathbf{y} \in \mathcal{D}_{II} \vee (\mathbf{y} \in \mathcal{D}_{III}, \: V_{CS}(\mathbf{y}) \geq t)), \: \theta_{CS}(\mathbf{y}) \leq t;\\
            \rho_{CC}(t, \mathbf{y}), \quad \text{otherwise}.
        \end{cases}
    \end{equation}
    Here
    \begin{equation}
        \rho_{CC}(t, \mathbf{y}) = \min_{\substack{\tau \in \{0, \tau_1(t, \mathbf{y}), \tau_2(t, \mathbf{y}), \tau_3(t, \mathbf{y})\} \\ 0 \leq \tau \leq \min(t, \frac{\pi}2)}} \sqrt{(\lvert x \rvert - x_{LR}(\tau, t))^2 + (y - y_{LR}(\tau, t))^2},
    \end{equation}
    \begin{equation}
        x_{LR}(\tau, t) \overset{\mathrm{def}}{=} 2\cos\tau - \cos(t - 2\tau) - 1, \quad y_{LR}(\tau, t) \overset{\mathrm{def}}{=} 2\sin\tau + \sin(t - 2\tau),
    \end{equation}
    \begin{equation}
        \tau_i(t, \mathbf{y}) \overset{\mathrm{def}}{=} \left(\frac{t}3 - 2\arctan\xi_i(t, \mathbf{y})\right) \: \mathrm{mod} \: 2\pi, \quad i \in \{1, 2, 3\},
    \end{equation}
    and $\xi_1(t, \mathbf{y})$, $\xi_2(t, \mathbf{y})$, $\xi_3(t, \mathbf{y})$ are real solutions of the cubic equation
    \begin{multline*}
        -\left(y + \sin\frac{t}3\right)\xi^3 + \left(3 + 3\lvert x \rvert + \cos\frac{t}3\right)\xi^2\\
        + \left(3y - \sin\frac{t}3\right)\xi + \cos\frac{t}3 - (1 + \lvert x \rvert) = 0.
    \end{multline*}
\end{assumption}
\begin{proof}
    It follows from the parametric description of $\mathcal{B}_{CS}(t)$, $\mathcal{B}_{CC}(t)$ from~\cite{Buzikov2021-md}. We left it to the reader to verify the solution of the real-valued function extreme problem.\qed
\end{proof}

\begin{lemma}\label{lem:tau_opt_dub}
    Let $t \in \mathbb{R}^+_0$, $\mathbf{y} \in \mathcal{Y} \setminus \mathcal{R}(t)$, $\rho(t, \mathbf{y}) > \ell$, and $\theta_{CS}(\mathbf{y}) \leq t$. If $\mathbf{y} \in \mathcal{D}_{II}$ or $\mathbf{y} \in \mathcal{D}_{III}$ and $V_{CS}(\mathbf{y}) \geq t$, then
    \begin{equation}
        T(t, \mathbf{y}) = \tau(t, \mathbf{y}) = t + \frac{V_{CS}(\mathbf{y}) - t - \ell}{1 + v}.
    \end{equation}
\end{lemma}
\begin{proof}
    Using the symmetry of the problem, we give the proof only for the case when $\mathbf{y}$ is on the right half plane. According to Eq. (3.2) from~\cite{Buzikov2021-md}, the point
    \begin{equation}
        \boldsymbol{\eta} =
        \begin{bmatrix}
            (t - \theta_{CS}(\mathbf{y}))\sin\theta_{CS}(\mathbf{y}) - \cos\theta_{CS}(\mathbf{y}) + 1\\
            (t - \theta_{CS}(\mathbf{y}))\cos\theta_{CS}(\mathbf{y}) + \sin\theta_{CS}(\mathbf{y})
        \end{bmatrix}
    \end{equation}
    can be reached by $CS$-trajectory with the velocity vector
    \begin{equation}
        \begin{bmatrix} \sin\theta_{CS}(\mathbf{y}) & \cos\theta_{CS}(\mathbf{y}) \end{bmatrix}^\top.
    \end{equation}
    This vector is co-directional with $\mathbf{y} - \boldsymbol{\eta}$, since
    \begin{equation}
        \mathbf{y} =
        \begin{bmatrix}
            (V_{CS}(\mathbf{y}) - \theta_{CS}(\mathbf{y}))\sin\theta_{CS}(\mathbf{y}) - \cos\theta_{CS}(\mathbf{y}) + 1\\
            (V_{CS}(\mathbf{y}) - \theta_{CS}(\mathbf{y}))\cos\theta_{CS}(\mathbf{y}) + \sin\theta_{CS}(\mathbf{y})
        \end{bmatrix}.
    \end{equation}
    The Dubins car moves with the maximal speed $\lVert \dot{\boldsymbol{y}}(t) \rVert = 1$. Thus, the target can be intercepted minimally in $(\lVert \mathbf{y} - \boldsymbol{\eta} \rVert - \ell) / (1 + v)$.\qed
\end{proof}

According to Lemma~\ref{lem:tau_is_conv}, sequence $\{t_n\}_{n = 0}^\infty$ given by $t_0 = 0$ and $t_n = \tau(t_{n - 1}, \boldsymbol{y}_T(t_{n - 1}))$ converges to $T^*[\boldsymbol{y}_T]$ for any $\boldsymbol{y}_T \in \mathrm{Lip}_v(\mathbb{R}^+_0, \mathcal{Y})$. Moreover, the step size $t_n - t_{n - 1}$ is maximal if $t = t_{n - 1}$ and $\mathbf{y} = \boldsymbol{y}_T(t_{n - 1})$ satisfy the requirements of Lemma~\ref{lem:tau_opt_dub}.

Comparing the obtained results with the known proposed algorithms~\cite{Clements1990-fp,Looker2008-yx,Zhang2014-fp,Meyer2015-st,Zheng2021-wu} for solving the interception by a Dubins car problem, we can conclude that Algorithm~\ref{alg:comp} works with a wider class of target movements (Lipschitz trajectories) and also guarantees convergence to the solution, if it exists.

\subsection{Numerical Experiments}\label{sec:numerical_experiments}

In this section, we present numerical experiments that solve the problem of the minimum-time interception of a moving target. The number of iterations required to reach the given precision for the optimal interception time is presented in Table~\ref{tab:comparison}. The rows in the table show an example of the target trajectory $\boldsymbol{y}_T \in \mathrm{Lip}_v(\mathbb{R}^+_0, \mathcal{Y})$. The desired precision $\delta \in \{10^{-3}, 10^{-6}, 10^{-9}\}$ is placed in the columns. The cells in the table contain the number of iterations $n$ for a given target trajectory and the plant model to reach the desired precision $T^*[\boldsymbol{y}_T] - t_n < \delta$.  It can be observed that the better precision, the larger number of iterations. As rule, if the target moves toward the reachable set, then the number of iterations is small. In contrast, if the target moves away from the reachable set, then the number of iterations is greater, than the speed of target greater.

\begin{table}
    \begin{center}
        \begin{minipage}{\textwidth}
            \caption{The number of iterations that have to be done to reach the given precision}\label{tab:comparison}
            \begin{tabular}{ccccccccc}
                \toprule
                & $\ell = \frac{1}{10}$ & & \multicolumn{3}{@{}c@{}}{Simple motions} & \multicolumn{3}{@{}c@{}}{Dubins car} \\
                \cmidrule{4-6}\cmidrule{7-9}
                & & $v$ & $10^{-3}$ & $10^{-6}$ & $10^{-9}$ &$10^{-3}$ & $10^{-6}$ & $10^{-9}$ \\
                \midrule
                \multirow{12}{*}{\rotatebox[origin=c]{90}{$\boldsymbol{y}_T(t) = \begin{bmatrix} \xi + vt\cos\varphi \\ \eta + vt\sin\varphi\end{bmatrix}$}} & \multirow{3}{0.12\linewidth}{$\xi = 0$ $\eta = 1$ $\varphi = 0$} & $1/4$ & $5$ &$10$ &$15$ &$5$ &$10$ &$15$ \\
                 & & $1/2$ &$10$ &$19$ &$29$ & $11$ &$23$ &$34$ \\
                 & & $3/4$ &$22$ &$43$ &$65$ & $48$ &$93$ & $137$ \\
                \cline{2-9}
                 & \multirow{3}{0.12\linewidth}{$\xi = 1$ $\eta = 1$ $\varphi = \pi/2$} & $1/4$ & $8$ &$15$ &$21$ &$7$ &$14$ &$20$ \\
                 & & $1/2$ &$17$ &$33$ &$48$ & $17$ &$32$ &$47$ \\
                 & & $3/4$ &$49$ &$90$ & $131$ & $49$ &$89$ & $130$ \\
                \cline{2-9}
                 & \multirow{3}{0.12\linewidth}{$\xi = -1$ $\eta = -2$ $\varphi = \pi/4$} & $1/2$ & $3$ & $5$ & $7$ & $11$ &$18$ &$25$ \\
                 & & $3/4$ & $3$ & $5$ & $8$ & $12$ &$20$ &$28$ \\
                 & & $1$ & $3$ & $6$ & $9$ & $14$ &$23$ &$33$ \\
                \cline{2-9}
                 & \multirow{3}{0.12\linewidth}{$\xi = -2$ $\eta = 0$ $\varphi = \pi/4$} & $1/2$ & $5$ & $9$ &$13$ & $19$ &$25$ &$30$ \\
                 & & $3/4$ & $6$ &$12$ &$18$ & $12$ &$31$ &$51$ \\
                 & & $1$ & $9$ &$18$ &$27$ &$5$ &$10$ &$15$ \\
                \cline{1-9}
                \cline{2-9}
                \multirow{16}{*}{\rotatebox[origin=c]{90}{$\boldsymbol{y}_T(t) = \begin{bmatrix} \xi + v/\omega_x\sin\omega_xt \\ \eta + v/\omega_y\sin\omega_yt\end{bmatrix}$}} & \multirow{4}{0.12\linewidth}{$\xi = 1$ $\eta = 1$ $\omega_x = 1$ $\omega_y = \sqrt{2}$} & $1/2$ & $5$ & $8$ &$11$ &$6$ & $9$ &$13$ \\
                 & & $1$ & $5$ & $7$ & $9$ &$7$ &$11$ &$16$ \\
                 & & $3/2$ & $5$ & $7$ & $8$ & $10$ &$20$ &$29$ \\
                 & & $2$ & $5$ & $7$ & $9$ & $28$ &$46$ &$64$ \\
                \cline{2-9}
                 & \multirow{4}{0.12\linewidth}{$\xi = -1$ $\eta = -2$ $\omega_x = 1$ $\omega_y = \sqrt{2}$} & $1/2$ &$12$ &$26$ &$40$ &$5$ & $8$ &$11$ \\
                 & & $1$ & $9$ &$21$ &$33$ &$6$ & $8$ &$10$ \\
                 & & $3/2$ & $7$ &$17$ &$26$ &$7$ &$10$ &$12$ \\
                 & & $2$ & $8$ &$20$ &$33$ &$9$ &$13$ &$16$ \\
                \cline{2-9}
                 & \multirow{4}{0.12\linewidth}{$\xi = -1$ $\eta = -2$ $\omega_x = 1$ $\omega_y = 2$} & $1/2$ &$11$ &$21$ &$30$ & $13$ &$23$ &$32$ \\
                 & & $1$ &$16$ &$26$ &$36$ & $19$ &$29$ &$39$ \\
                 & & $3/2$ &$18$ &$26$ &$33$ & $21$ &$28$ &$36$ \\
                 & & $2$ &$20$ &$26$ &$31$ & $25$ &$31$ &$36$ \\
                \cline{2-9}
                 & \multirow{4}{0.12\linewidth}{$\xi = 0$ $\eta = -1$ $\omega_x = 2$ $\omega_y = 1$} & $1/2$ & $3$ & $6$ & $9$ &$9$ &$14$ &$18$ \\
                 & & $1$ & $6$ &$12$ &$19$ & $12$ &$17$ &$22$ \\
                 & & $3/2$ &$17$ &$37$ &$57$ & $17$ &$22$ &$27$ \\
                 & & $2$ &$21$ &$36$ &$51$ &$9$ &$16$ &$23$ \\
                 \bottomrule
            \end{tabular} 
        \end{minipage}
    \end{center}
\end{table}

\begin{figure}
    \centering
    \includegraphics[width=0.59\textwidth]{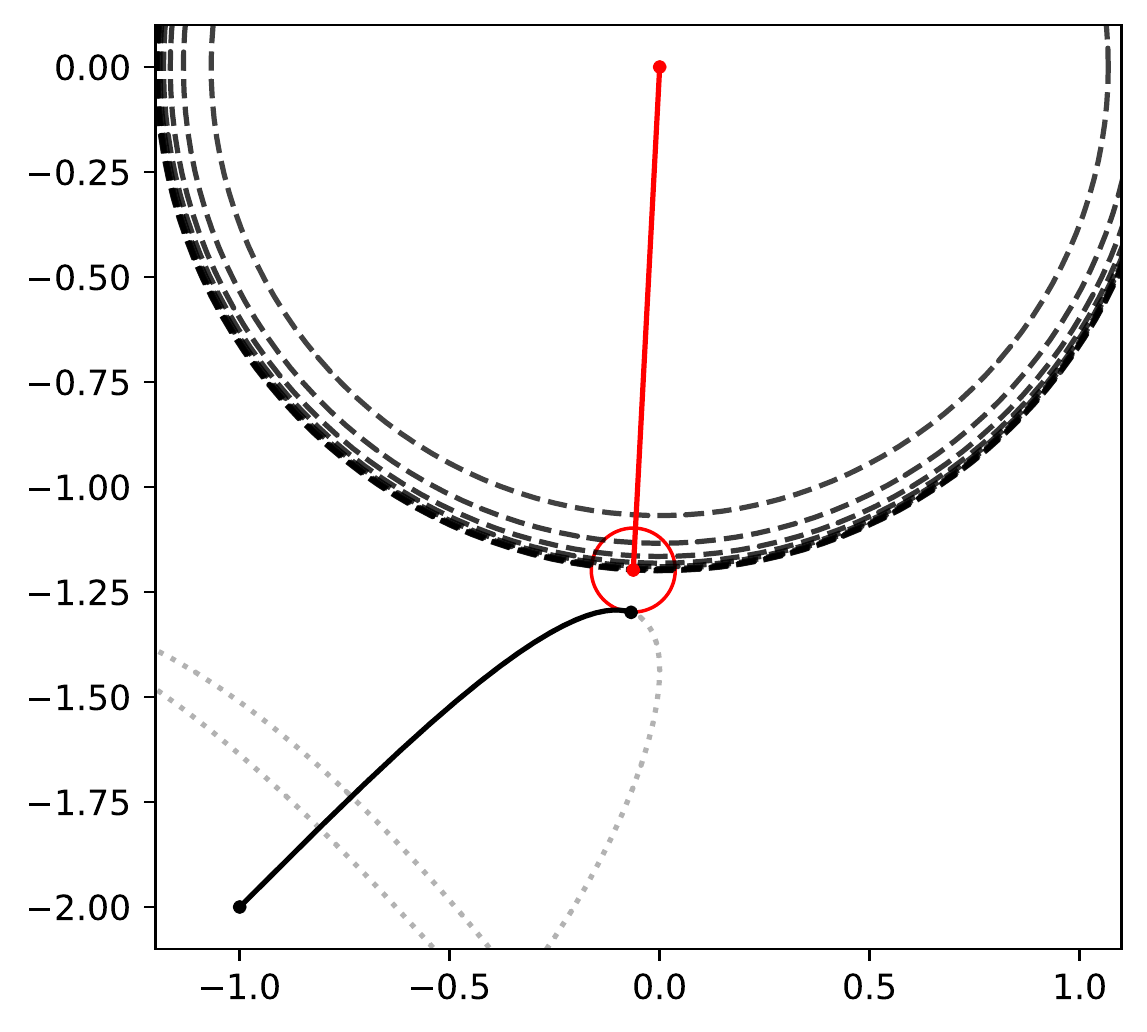}
    \includegraphics[width=0.40\textwidth]{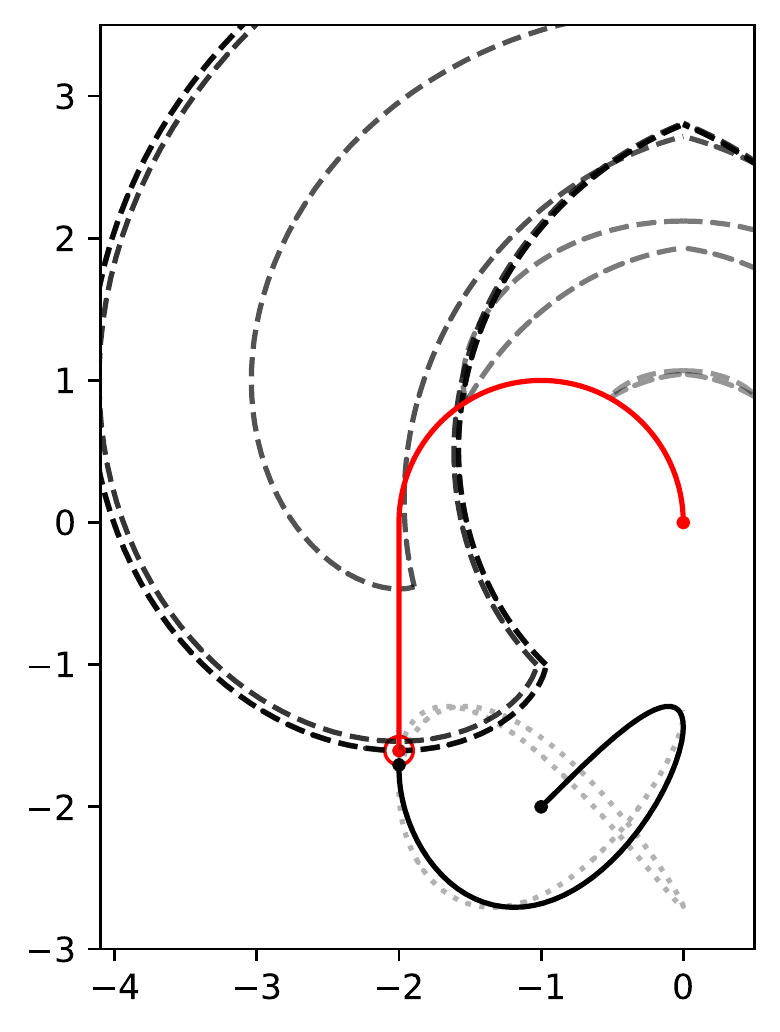}
    \caption{The paths of time-minimal interception (red lines) of the target that moves along the Lissajous curve (black solid and dotted lines) described in Table~\ref{tab:comparison} ($\xi = -1$, $\eta = -2$, $\omega_x = 1$, $\omega_y = \sqrt{2}$, $v = 1$). The dashed lines correspond to the planar reachable set at the time moments produced by the fixed-point algorithm for the simple motions (left) and the Dubins car (right). The red circles represent the capture sets $\ell = 1/10$. According to Table~\ref{tab:comparison}, nine iterations were required to reach a precision $10^{-3}$ in the simple motion case and six iterations were required in the Dubins car case.}
    \label{fig:path_ex}
\end{figure}

Examples of optimal paths that intercept the moving target for a given trajectory $\boldsymbol{y}_T$ are presented in Fig~\ref{fig:path_ex}.

\section{Conclusion and Further Work}\label{sec:conclusion}

We presented a study of the minimum-time interception problem by employing a reachable set analysis. We showed that if the distance from an arbitrary point to the reachable set projection is a given function, then the minimal-time interception can be calculated as the smallest non-negative root of the equation in which the distance from the reachable set projection to the position of the moving target is equal to the capture radius (Theorem~\ref{th:main_equation}). We propose two functions (Definitions~\ref{def:best_est} and~\ref{def:tau}) that can be used in a fixed-point iteration method, which always converges to the optimal solution (Lemma~\ref{lem:tau_is_conv} and Theorem~\ref{th:T_corr_conv}). The function given by Definition~\ref{def:tau} can be explicitly calculated in terms of the distance to the reachable set projection. In turn, the function given by Definition~\ref{def:best_est} can be calculated as a root of a real equation (Theorem~\ref{th:T_calc}). In practice, this equation can be more difficult to solve than the original problem. However, if there is a way to effectively calculate a solution to this equation, then the step of the fixed-point iteration method based on the function from Definition~\ref{def:best_est} is the largest, which can be performed for an arbitrary Lipschitz target trajectory without losing the guarantee of not missing the least root.

For the Dubins model, we obtained an explicit analytical form of the distance to the planar reachable set (Assumption~\ref{assum:rho}). In addition, we proved that the simple universal lower estimator is equal to the best universal lower estimator in a significant area of the state space (Lemma~\ref{lem:tau_opt_dub}).

An interesting extension of the presented results would be the generalization of the framework to metric spaces. The current results operate with finite-dimensional normed spaces, but some practically important problems, such as {\it lateral} interception of a moving target by a Dubins car, cannot be properly formalized in the normed spaces. In addition, the isotropic rocket model seems perspective to imply all aspects of the proposed algorithmic framework.

~

\textbf{Data availability} All data generated or analysed during this study are included in this published article.

~

\textbf{Conflict of interest} The author has no conflicts of interest to declare.

\begin{acknowledgements}
    The research was supported by RSF (project No. 23-19-00134).
\end{acknowledgements}

\bibliographystyle{spmpsci}
\bibliography{main}

\end{document}